# THE DOUBLE-CONSTANT MATRIX, CENTERING MATRIX AND EQUICORRELATION MATRIX: THEORY AND APPLICATIONS

B. O'Neill,* *Australian National University*\*\*

WRITTEN 7 AUGUST 2019; REVISED 3 DECEMBER 2020

**Abstract**

This paper examines the properties of real symmetric square matrices with a constant value for the main diagonal elements and another constant value for all off-diagonal elements. This matrix form is a simple subclass of circulant matrices, which is a subclass of Toeplitz matrices. It encompasses other useful matrices such as the centering matrix and the equicorrelation matrix, which arise in statistical applications. We examine the general form of this class of matrices and derive its eigendecomposition and other important properties. We use this as a basis to look at the properties of the centering matrix and the equicorrelation matrix, and various statistics that use these matrices.

DOUBLE-CONSTANT MATRICES; CENTERING MATRIX; EQUICORRELATION MATRIX; CANONICAL FORM.

Matrix algebra forms an important part of statistical modelling and is a substantial part of the basis of linear models in statistics. So important is the subject that there are a number of books on matrix algebra specifically for the statistical profession (see e.g., Searle 2006, Harville 2008, Seber 2008, Gentle 2017). In this specialist context, in addition to setting out the general theory of matrices, special attention is often paid to particular matrix forms that commonly arise in statistical applications. The value of focussing on special matrix forms is that it allows a deep examination of properties and the derivation of useful formulae for analysis.

In this paper we will examine the properties of a simple class of real symmetric square matrices with a fixed constant on the main diagonal, and another fixed constant on all the off-diagonal elements (we will call these **double-constant matrices**). This is a simple matrix form which is a subclass of circulant matrices, which is itself a subclass of Toeplitz matrices. There is a well-developed mathematical literature on Toeplitz and circulant matrices (Grenander and Szegö 1958; Davis 1979; Böttcher and Silverman 1999; Böttcher and Grundsky 2000, 2005; Gray 2006), but in this paper we will examine this smaller and simpler subclass of matrices in its own right, with a view to giving insight into its properties and applications. While this matrix form is *prima facie* quite simple (some might even say trivial) a deep examination of its properties is warranted by the fact that it frequently appears in statistical applications, and its properties elucidate the behaviour of some important statistics in these applications.

---

* E-mail address: ben.oneill@hotmail.com.
\*\* Research School of Population Health, Australian National University, Canberra ACT 0200, Australia.

The double-constant matrix arises in statistical problems via two main types of matrices, which are the centering matrix and the equicorrelation matrix. Both these forms have a constant value down the main diagonal and another (non-zero) constant on the off-diagonal elements. We will see that these matrices form important cases of the class of double-constant matrices, and we will explain why these two particular matrix forms arise in statistical problems. We defer presentation of literature on these latter matrices until after we have examined the properties of the double-constant matrix more generally.

Statistical texts on matrix algebra give special attention to the centering matrix, and most also give some attention to the equicorrelation matrix. In some texts there are passing exercises or discussion that looks at the double-constant matrix, though the form is not examined in detail. For example, Seber (2008) devotes a short section to matrices with repeated elements and block (§15.4, pp. 316-319) and this part includes a brief discussion of the double-constant matrix.[1] In this paper we will go further, giving a detailed examination of the general form of the double-constant matrix, and showing a number of applications of this matrix in statistical problems. We will establish the eigendecomposition of the matrix form, and establish its properties under simple matrix operations and matrix functions. We will examine the possible values of the eigenvectors of the matrix and use this to decompose the class of double-constant matrices into important subgroups. Through this analysis we will see that the centering matrix and the equicorrelation matrix arise as important representatives of general classes of double-constant matrices, grouped by their definiteness properties.

In addition to examining the double-constant form in general, we will establish a number of useful results pertaining to the centering matrix and equicorrelation matrix. We will show the connection of the centering matrix to the concept of "degrees of freedom" in linear models and we will show some decomposition results that are useful for sum-of-squares and sample variance decompositions. Similarly, we will show the connection of the equicorrelation matrix to the notion of the "effective sample size" and "effective degrees of freedom" in models with equicorrelated data. We hope that this broad examination will shed some light on an important class of matrices that are simple in their basic form, but rich in statistical applications.

---

[1] Seber does not name this matrix form, and so he does not refer to it as the "double constant matrix" (that is our name for it). His discussion shows the general form of the matrix and gives some results for the matrix.

# 1. The double-constant matrix and its eigendecomposition

There are a several useful matrices in statistics that are square matrices with diagonal elements that are constant, and off-diagonal elements that are also constant (these two constants being able to differ from one another). These matrices can be formed as the sum of a square diagonal matrix and a square constant matrix, leading to an $n \times n$ real symmetric matrix:

$$\mathbf{M} \equiv \mathbf{M}(a, t) \equiv (a - t) \cdot \mathbf{I}_{n \times n} + t \cdot \mathbf{1}_{n \times n} = \begin{bmatrix} a & t & t & \cdots & t \\ t & a & t & \cdots & t \\ t & t & a & \cdots & t \\ \vdots & \vdots & \vdots & \ddots & \vdots \\ t & t & t & \cdots & a \end{bmatrix} \qquad a, t \in \mathbb{R}.$$

Despite its simple and useful structure, this matrix form does not appear to have s name in the mathematics and statistics literature. It is a special case of the broader form of the circulant matrices, which is a form of Toeplitz matrices (also called diagonal-constant matrices). The present matrix form is simpler, insofar as all the secondary diagonals have the same constant value. We will call this matrix form the **double-constant matrix**.

**REMARK:** Every double-constant matrix is a scalar multiple of the double-constant matrix with $a = 1$. Thus, it is possible to remove the scalar multiple and proceed with $a = 1$ without any serious loss of generality. Nevertheless, we will proceed for the above case where the scalar multiple is built into the matrix. This turns out to be more convenient for later analysis. □

We will derive a number of the important properties of this matrix, including its determinant, characteristic polynomial, eigenvalues and eigenvectors, and is eigendecomposition. This will allow us to take powers of the matrix, including its inverse and its various principal roots. It also allows us to obtain simple expressions for broader functions of the matrix. To facilitate our analysis, we note in advance that the matrix has two eigenvalues (ignoring multiplicity), which are the **major eigenvalue** $\lambda_* \equiv a - t$ and **minor eigenvalue** $\lambda_{**} \equiv a - t + nt$. (These eigenvalues are distinct so long as $t \neq 0$.) The element values in the double-constant matrix can be recovered from the eigenvalues as:

$$a = \frac{\lambda_{**} + (n-1)\lambda_*}{n} \qquad t = \frac{\lambda_{**} - \lambda_*}{n}.$$

Hence, we can write the double-constant matrix in its **canonical form** (we abuse notation by using the same notation for matrix functions using the constant values or eigenvalues) as:

$$\mathbf{M} \equiv \mathbf{M}(\lambda_*, \lambda_{**}) \equiv \lambda_* \cdot \left( \mathbf{I}_{n \times n} - \frac{\mathbf{1}_{n \times n}}{n} \right) + \lambda_{**} \cdot \frac{\mathbf{1}_{n \times n}}{n}.$$

**THEOREM 1:** For all $a, t \in \mathbb{R}$ the double-constant matrix has determinant:

$$\det(\mathbf{M}) = \lambda_*^{n-1} \lambda_{**} = (a-t)^{n-1}(a-t+nt).$$

**THEOREM 2:** For all $a, t \in \mathbb{R}$ the characteristic polynomial of the double-constant matrix is:

$$p(\lambda) \equiv \det(\mathbf{M} - \lambda \cdot \mathbf{I}_{n \times n}) = (\lambda_* - \lambda)^{n-1}(\lambda_{**} - \lambda) = (a - t - \lambda)^{n-1}(a - t + nt - \lambda).$$

The eigenvalues of the double-constant matrix are:

$$\lambda_1 = \lambda_{**} = a - t + nt \qquad \lambda_2 = \cdots = \lambda_n = \lambda_* = a - t.$$

In the case where these eigenvalues are non-zero, the corresponding eigenvectors are:

| Eigenvalues | Eigenvectors |
|---|---|
| $\lambda_{**}$ | Every vector $\mathbf{u}_{**} \propto \mathbf{1}_n$, |
| $\lambda_*$ | Every non-zero vector $\mathbf{u}_*$ with $\mathbf{u}_* \cdot \mathbf{1}_n = 0$. |

**THEOREM 3:** Define the $n \times n$ unitary discrete Fourier transform (DFT) matrix:

$$\mathbf{U} \equiv \frac{1}{\sqrt{n}} \begin{bmatrix} 1 & 1 & 1 & \cdots & 1 \\ 1 & \omega_n^1 & \omega_n^2 & \cdots & \omega_n^{n-1} \\ 1 & \omega_n^2 & \omega_n^4 & \cdots & \omega_n^{2(n-1)} \\ \vdots & \vdots & \vdots & \ddots & \vdots \\ 1 & \omega_n^{n-1} & \omega_n^{2(n-1)} & \cdots & \omega_n^{(n-1)(n-1)} \end{bmatrix} \qquad \omega_n \equiv \exp\left(-\frac{2\pi i}{n}\right).$$

(This is a symmetric unitary matrix so $\mathbf{U}^{-1} = \mathbf{U}^* = \overline{\mathbf{U}}^{\mathrm{T}} = \overline{\mathbf{U}}$, where the "bar" operation denotes the complex conjugate of the original matrix. The latter matrix is the inverse-DFT matrix.) The double-constant matrix is diagonalized by this matrix, and can be written as:

$$\mathbf{M} = \mathbf{U} \boldsymbol{\Lambda} \overline{\mathbf{U}} = \overline{\mathbf{U}} \boldsymbol{\Lambda} \mathbf{U} \qquad \boldsymbol{\Lambda} \equiv \mathrm{diag}(\lambda_{**}, \lambda_*, \ldots, \lambda_*).$$

We have now established that the double-constant matrix is a diagonalizable matrix that can be fully characterised by its major and minor eigenvalues (which are distinct if $t \neq 0$). The matrix is diagonalized by the unitary DFT matrix, which is actually a general property of all circulant matrices (see Gray 2006, § 3). The rows/columns of the DFT matrix and inverse-DFT matrix (its conjugate transpose) are eigenvectors of the double-constant matrix. The first row/column of the unitary DFT matrix is the unit eigenvector for the minor eigenvalue and the other rows/columns are orthonormal eigenvectors corresponding to the major eigenvalue (all multiplicities). This diagonalization establishes an alternative way of looking at the connection between the constants in the double-constant matrix and its eigenvalues. We can rearrange the diagonalization to obtain $\mathbf{U}\mathbf{M} = (\mathbf{U}\boldsymbol{\Lambda})^{\mathrm{T}}$ and $\overline{\mathbf{U}}\mathbf{M} = (\overline{\mathbf{U}}\boldsymbol{\Lambda})^{\mathrm{T}}$, so the DFT (inverse-DFT) of the double-constant matrix is the transpose of the DFT (inverse-DFT) of its eigenvalue matrix.

Like any square matrix, the double-constant matrix can be treated as a linear transformation $\mathbf{M}: \mathbb{R}^{n \times m} \to \mathbb{R}^{n \times m}$ operating on a matrix $\pmb{x} \in \mathbb{R}^{n \times m}$.[2] This transformation occurs via matrix pre-multiplication $\pmb{x} \mapsto \mathbf{M}\pmb{x}$. The eigendecomposition of $\mathbf{M}$ uses the unitary DFT matrix $\mathbf{U}$ as its eigenvector matrix, so there is a close connection between transformation using the double-constant matrix, and corresponding operations conducted in Fourier space. To understand the double-constant matrix in greater detail, it is worth examining these corresponding operations.

Consider an arbitrary data vector $\pmb{x}_k$ (composed of $n$ elements) which is a column of the matrix $\pmb{x}$. Let $\mathcal{F}_{\pmb{x}_k}(r)$ denote the unitary DFT of $\pmb{x}_k$ taken at the frequency $2\pi r/n$.[3] Then taking the DFT of the entire matrix $\pmb{x} = [\pmb{x}_1 \cdots \pmb{x}_m]$ we obtain:

$$\mathbf{U}\pmb{x} = \begin{bmatrix} \mathcal{F}_{\pmb{x}_1}(0) & \mathcal{F}_{\pmb{x}_2}(0) & \cdots & \mathcal{F}_{\pmb{x}_m}(0) \\ \mathcal{F}_{\pmb{x}_1}(1) & \mathcal{F}_{\pmb{x}_2}(1) & \cdots & \mathcal{F}_{\pmb{x}_m}(1) \\ \vdots & \vdots & \ddots & \vdots \\ \mathcal{F}_{\pmb{x}_1}(n-1) & \mathcal{F}_{\pmb{x}_2}(n-1) & \cdots & \mathcal{F}_{\pmb{x}_m}(n-1) \end{bmatrix}.$$

Since $\mathbf{\Lambda} = \text{diag}(\lambda_{**}, \lambda_*, \ldots, \lambda_*)$ we then have:

$$\mathbf{\Lambda}\mathbf{U}\pmb{x} = \begin{bmatrix} \lambda_{**} \cdot \mathcal{F}_{\pmb{x}_1}(0) & \lambda_{**} \cdot \mathcal{F}_{\pmb{x}_2}(0) & \cdots & \lambda_{**} \cdot \mathcal{F}_{\pmb{x}_m}(0) \\ \lambda_* \cdot \mathcal{F}_{\pmb{x}_1}(1) & \lambda_* \cdot \mathcal{F}_{\pmb{x}_2}(1) & \cdots & \lambda_* \cdot \mathcal{F}_{\pmb{x}_m}(1) \\ \vdots & \vdots & \ddots & \vdots \\ \lambda_* \cdot \mathcal{F}_{\pmb{x}_1}(n-1) & \lambda_* \cdot \mathcal{F}_{\pmb{x}_2}(n-1) & \cdots & \lambda_* \cdot \mathcal{F}_{\pmb{x}_m}(n-1) \end{bmatrix}.$$

Finally, we have the transformed value $\mathbf{M}(\pmb{x}) = \mathbf{M}\pmb{x} = \overline{\mathbf{U}}\mathbf{\Lambda}\mathbf{U}\pmb{x} = \overline{\mathbf{U}}(\mathbf{\Lambda}\mathbf{U}\pmb{x})$, which means that we transform $\pmb{x}$ by taking the inverse-DFT of the above matrix. This show us that applying the double-constant matrix as a linear transform to a matrix $\pmb{x}$ is equivalent to scaling the Fourier values of the columns of $\pmb{x}$. The Fourier values at the zero frequency are scaled by the minor eigenvalue $\lambda_{**}$, and the Fourier values at all the non-zero frequencies are scaled by the major eigenvalue $\lambda_*$. Intuitively, this result is unsurprising, in view of the original form of the double-constant matrix. The double-constant matrix is a linear combination of the identity matrix and the unit matrix, so it transforms a vector by scaling its original values (using the identity matrix) and scaling the sum of its values (using the unit matrix), and adding these two parts together. In Fourier space the first of these operations corresponds to scaling all the frequencies, and the latter corresponds to scaling only the zero frequency.

---

[2] Actually, the matrix can just as easily be treated as a linear transform operating on *complex* matrices, but we will confine our attention in this paper to real matrices, in view of our focus on statistical data.

[3] To avoid ambiguity, if $\pmb{x}_k$ is composed of elements $x_{k,1}, \ldots, x_{k,n}$ we take $\mathcal{F}_{\pmb{x}_k}: \mathbb{Z} \to \mathbb{C}$ to be defined by:

$$\mathcal{F}_{\pmb{x}_k}(r) \equiv \frac{1}{\sqrt{n}} \sum_{t=0}^{n-1} x_{k,t+1} \cdot \exp\left(-\frac{2\pi i r t}{n}\right) \qquad \text{for all } r \in \mathbb{Z}.$$

In the above analysis we have looked at linear transformation via pre-multiplication $x \mapsto \mathbf{M}x$, but the same basic dualism emerges when we look at post-multiplication $x \mapsto x\mathbf{M}$. (This leads to a similar form that can easily be derived by analogy to the above method; this is left as an exercise by the reader.) In either case, we obtain a dual interpretation of the double-constant matrix as a linear transform. It can be viewed as a transformation operating in Euclidean space as a linear combination of the identity matrix and unit matrix, or it can be viewed equivalently as a transformation operating in Fourier space as a combination of scaling of the signal power at the zero frequency and separate scaling of the signal power at the non-zero frequencies.

It is worth noting that this dualism between the Euclidean and Fourier space allows us to appeal to all the various properties of discrete Fourier transformation, including standard theorems like Parseval's theorem and the Plancherel theorem, that relate the dot products of transformed vectors in Euclidean space to conjugated products of the signal in Fourier space. In particular, for any vectors $x, y \in \mathbb{R}^n$ that are subjected to the double-constant transformations $\mathbf{M}_1$ and $\mathbf{M}_2$ respectively, we have this manifestation of Parseval's theorem:

$$\begin{aligned}
\mathbf{M}_1(x) \cdot \mathbf{M}_2(y) &= [\mathbf{M}_1 x]^\mathrm{T} [\mathbf{M}_2 y] \\
&= [\bar{\mathbf{U}} \Lambda_1 \mathbf{U} x]^\mathrm{T} [\mathbf{U} \Lambda_2 \bar{\mathbf{U}} y] \\
&= x^\mathrm{T} \mathbf{U} \Lambda_1 \bar{\mathbf{U}} \mathbf{U} \Lambda_2 \bar{\mathbf{U}} y \\
&= x^\mathrm{T} \mathbf{U} \Lambda_1 \Lambda_2 \bar{\mathbf{U}} y \\
&= (\mathbf{U} x)^\mathrm{T} \Lambda_1 \Lambda_2 (\bar{\mathbf{U}} y) \\
&= \begin{bmatrix} \mathcal{F}_x(0) \\ \mathcal{F}_x(1) \\ \vdots \\ \mathcal{F}_x(n-1) \end{bmatrix}^\mathrm{T} \begin{bmatrix} \lambda_{1**}\lambda_{2**} & 0 & \cdots & 0 \\ 0 & \lambda_{1*}\lambda_{2*} & \cdots & 0 \\ \vdots & \vdots & \ddots & \vdots \\ 0 & 0 & \cdots & \lambda_{1*}\lambda_{2*} \end{bmatrix} \begin{bmatrix} \overline{\mathcal{F}_y(0)} \\ \overline{\mathcal{F}_y(1)} \\ \vdots \\ \overline{\mathcal{F}_y(n-1)} \end{bmatrix} \\
&= \lambda_{1**}\lambda_{2**} \cdot \mathcal{F}_x(0) \cdot \overline{\mathcal{F}_y(0)} + \lambda_{1*}\lambda_{2*} \sum_{t=1}^{n-1} \mathcal{F}_x(t) \cdot \overline{\mathcal{F}_y(t)},
\end{aligned}$$

and this manifestation of the Plancherel theorem:

$$\begin{aligned}
\|\mathbf{M}_1(x)\|^2 &= \mathbf{M}_1(x) \cdot \mathbf{M}_1(x) \\
&= \lambda_{1**}^2 \cdot \mathcal{F}_x(0) \cdot \overline{\mathcal{F}_x(0)} + \lambda_{1*}^2 \sum_{t=1}^{n-1} \mathcal{F}_x(t) \cdot \overline{\mathcal{F}_x(t)} \\
&= \lambda_{1**}^2 \cdot \|\mathcal{F}_x(0)\|^2 + \lambda_{1*}^2 \sum_{t=1}^{n-1} \|\mathcal{F}_x(t)\|^2 .
\end{aligned}$$

## 2. Further properties and closure under operations

The dualistic interpretation of transformation via the double-constant matrix means that there are two simple ways to deal with operations involving this matrix form. One is to substitute the expansion of the double-constant matrix as a linear combination of the identity matrix and the unit matrix, and perform operations on these latter matrices. Another is to appeal to the eigendecomposition of the double-constant matrix and take advantage of the fact that all such matrices have a common orthonormal eigenvector matrix given by the DFT matrix and its inverse.[4] Either of these techniques can be used to derive relevant results for operations on this type of matrix; the first method deals directly with the constants in the matrix, while the second transitions through their eigenvalues. In order to illustrate this duality, and some useful results from using each method, we present two general rules for sums and products of double-constant matrices, derived via each of these respective methods.

**THEOREM 4:** For all constants $\kappa_i, a_i, t_i$ for $i = 1, \ldots, m$ we have the linear combination:

$$\sum_{i=1}^{m} \kappa_i \cdot \mathbf{M}(a_i, t_i) = \mathbf{M}\left(\sum_{i=1}^{m} \kappa_i \cdot a_i, \sum_{i=1}^{m} \kappa_i \cdot t_i\right).$$

**THEOREM 5:** For all constants $a_i, t_i$ for $i = 1, \ldots, m$ we have the product:

$$\prod_{i=1}^{m} \mathbf{M}(a_i, t_i) = \mathbf{M}(a_\times, t_\times),$$

where the new constants in the product-matrix are:

$$a_\times = \frac{1}{n}\left(\prod_{i=1}^{m}(a_i - t_i + nt_i) + (n-1)\prod_{i=1}^{m}(a_i - t_i)\right),$$

$$t_\times = \frac{1}{n}\left(\prod_{i=1}^{m}(a_i - t_i + nt_i) - \prod_{i=1}^{m}(a_i - t_i)\right).$$

(This result shows that products of double-constant matrices are commutative and associative. We can therefore refer to arbitrary products unambiguously using the product notation used for scalar products, just as we have done here.)

---

[4] The unitary DFT matrix does not depend on the values $a$ or $t$, which means that this matrix is the orthonormal eigenvector matrix for *all possible double-constant matrices*. (In fact, this is true not just for double-constant matrices, but for all circulant matrices.)

As with our other results, proofs of these theorems are shown in the Appendix. However, we note that the first result is easily derived from the decomposition of the double-constant matrix in Euclidean space, and the second result is easily derived by transformation of the eigenvectors in Fourier space. Either operation is handled easily in Fourier space (i.e., by operations on the eigenvalues), but the linear combinations are also handled easily in Euclidian space (i.e., by operations on the constants in the matrix). It is quite simple to derive further results for linear combinations of powers of double-constant matrices (i.e., polynomials of these matrices), and other combinations of double-constant matrices. This matrix class is amenable to many simple operations since it is closed under a range of useful operation types, including summation, multiplication, powers, exponentials, inversion, square-roots, etc.

Another important consequence of the eigendecomposition of the double-constant matrix (which follows from Sylvester's theorem) is that every analytic function of a double-constant matrix yields another double-constant matrix as its output. This occurs because eigenvectors for the double-constant matrix are fixed across all matrices of this kind; analytic functions change the eigenvalues in the decomposition, but leave the eigenvectors unchanged, which means that we still have a double-constant matrix. Since the major and minor eigenvalue map directly to the constants in the matrix, this means that analytic functions change the constants in the double-constant matrix, but do not change the matrix form. This makes it easy to deal with the double-constant matrix in problems involving analytic transformations.

**THEOREM 6:** Consider an arbitrary analytic function $f$ where $f(\lambda_{**})$ and $f(\lambda_*)$ are well-defined. Then we have $f(\mathbf{M}(a,t)) = \mathbf{M}(\tilde{a}, \tilde{t})$ where the new constants for the matrix are:

$$\tilde{a} \equiv \frac{f(\lambda_{**}) + (n-1)f(\lambda_*)}{n} = \frac{f(a - t + nt) + (n-1)f(a - t)}{n},$$

$$\tilde{t} \equiv \frac{f(\lambda_{**}) - f(\lambda_*)}{n} = \frac{f(a - t + nt) - f(a - t)}{n}.$$

We can write this transformed matrix explicitly as:

$$f(\mathbf{M}(a,t)) = f(a - t) \cdot \mathbf{I}_{n \times n} + \frac{f(a - t + nt) - f(a - t)}{n} \cdot \mathbf{1}_{n \times n}.$$

We can see from this theorem that the class of double-constant matrices is closed under analytic transformations. When we apply an analytic transform to the double-constant matrix, although the constant values in the matrix change (according to the above formulae), the basic form of the matrix remains the same. Some particular functions of interest are the exponential of the matrix and the logarithm of the matrix, which are given respectively by:

$$e^{\mathbf{M}} = e^{a-t}\left[\mathbf{I}_{n\times n} + \frac{1}{n}(e^{nt} - 1)\cdot \mathbf{1}_{n\times n}\right],$$

$$\ln \mathbf{M} = \ln(a-t)\left[\mathbf{I}_{n\times n} + \frac{1}{n}\left(\frac{\ln(a-t+nt)}{\ln(a-t)} - 1\right)\cdot \mathbf{1}_{n\times n}\right].$$

(The latter exists only if the double-constant matrix is non-singular. If $\lambda_* > 0$ and $\lambda_{**} > 0$ then there is a unique real principal logarithm.) The powers of the matrix are given by:

$$\mathbf{M}^\gamma = (a-t)^\gamma\left[\mathbf{I}_{n\times n} + \frac{1}{n}\left(\left(1+\frac{nt}{a-t}\right)^\gamma - 1\right)\cdot \mathbf{1}_{n\times n}\right].$$

If $\lambda_* \ne 0$ and $\lambda_{**} \ne 0$ then the double-constant matrix is invertible, with inverse:

$$\mathbf{M}^{-1} = \frac{1}{a-t}\left[\mathbf{I}_{n\times n} - \frac{t}{a-t+nt}\cdot \mathbf{1}_{n\times n}\right].$$

Another power of interest is the principal square-root of the double-constant matrix, which is:

$$\mathbf{M}^{1/2} = \sqrt{a-t}\left[\mathbf{I}_{n\times n} + \frac{1}{n}\left(\sqrt{\frac{a-t+nt}{a-t}} - 1\right)\cdot \mathbf{1}_{n\times n}\right].$$

The reader will have noticed that we have not imposed limitations on the constants $a$ and $t$, though some of the results —such as the matrix inverse— require restrictions. This means that the "double-constant" matrix, as presently specified, encompasses the case $a = t$ where we really have only a constant in the matrix, and the case $a \ne t = 0$ where we have a scaled identity matrix. If we wish to exclude these cases from the ambit of discussion, we will refer to the double-constant matrix with $a \ne t$ as a **proper double-constant matrix** and we refer to the double-constant matrix with $a \ne t \ne 0$ as a **non-degenerate double-constant matrix**.

It is easy to see that every proper double-constant matrix has a non-zero major eigenvalue, and we can further demarcate the class of non-degenerate double-constant matrices by grouping them into the singular and non-singular cases, corresponding to whether or not the minor eigenvalue is zero. In the case where the minor eigenvalue is zero, the non-degenerate double-constant matrix is proportional to the centering matrix $\mathbf{C}_{n\times n} = \mathbf{I}_{n\times n} - \mathbf{1}_{n\times n}/n$ (which is itself a particular case of the double-constant matrix), and in the case where the minor eigenvalue is non-zero we can further partition the matrix class into negative definite, indefinite, and positive definite forms, the latter giving matrices that are proportional to the equicorrelation matrix used in statistics. In Figure 1 below we show this partition of the matrix class.

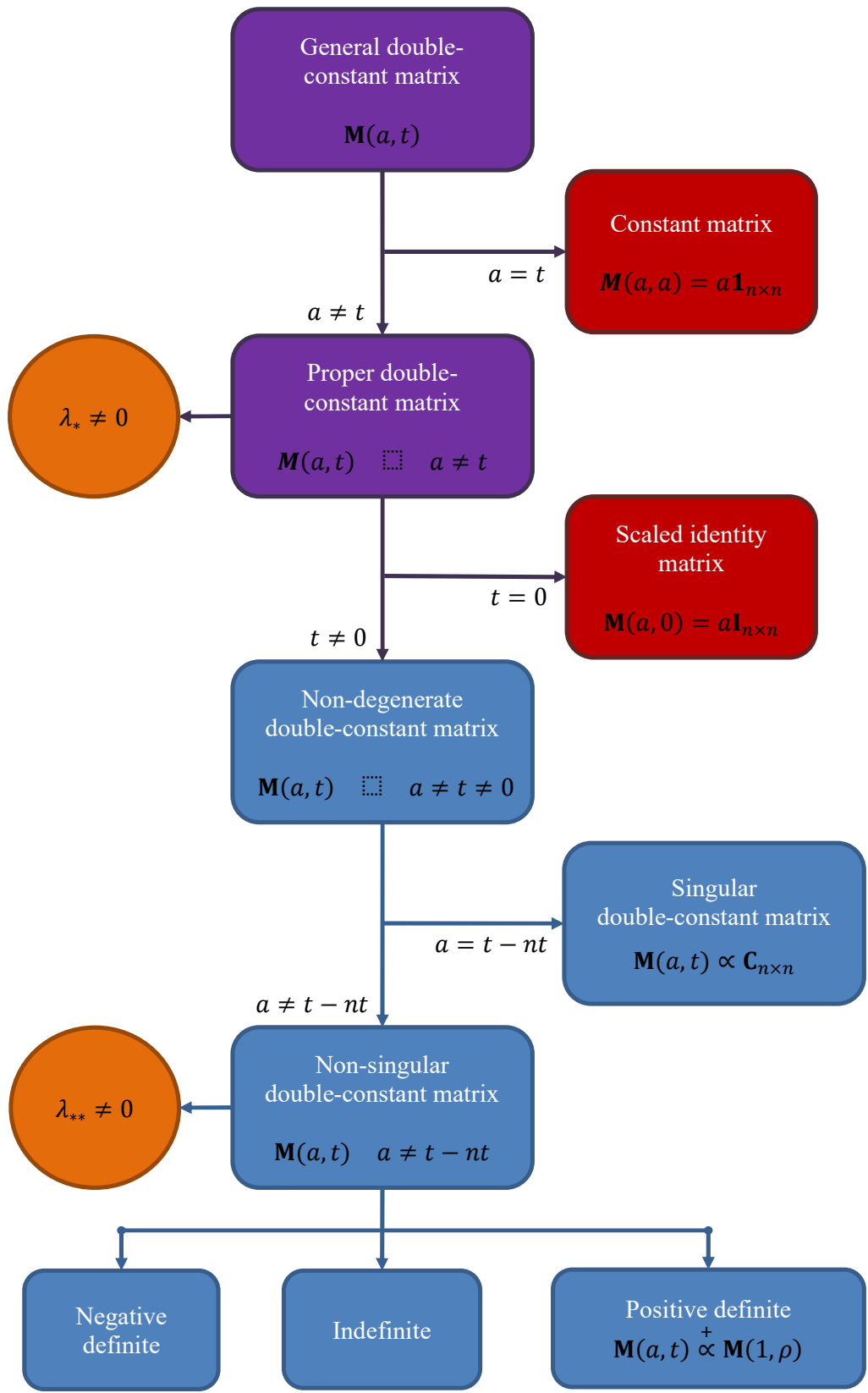

**FIGURE 1** – Partitioning the class of double-constant matrices

As can be seen from the figure, the partition we have imposed here is based on partitioning the cases where the major and minor eigenvalues are zero, negative, or positive, with some other exclusions for degenerate cases. We first narrow the class of double-constant matrices by excluding constant matrices (where the major eigenvalue is zero) and scaled identity matrices (where the eigenvalues are identical), leaving us with the subclass of non-degenerate double-constant matrices. It is easy to see that every non-degenerate double-constant matrix has full rank, though lower rank occurs when we have either a non-zero constant matrix (with unit rank) or the zero matrix (with zero rank).

Once we have narrowed down to the non-degenerate double-constant matrix we then partition into the singular case where we have a matrix proportional to the centering matrix, and the non-singular case, which includes negative definite, indefinite, and positive definite classes. The two parts of this partition that are of interest in statistical applications are the singular non-degenerate double-constant matrix (proportional to the **centering matrix**) and the positive definite non-degenerate double-constant matrix (positively proportional to the **equicorrelation matrix**). Both of these matrices are special cases of the class of double-constant matrices, but they also represent archetypal elements of these partition subclasses, insofar as all matrices in those subclasses are merely multiples of these archetypal matrices.

## 3. The centering matrix and its uses in statistical models

In statistical models it is common to wish to center columns of values in a matrix (or a column vector) by shifting their location to have zero mean. This is done in the context of regression models, where much of the mathematical analysis is simplified by centering the explanatory vectors prior to analysis. Centering of data is done by projecting a matrix onto the column space orthogonal to the unit vector, which is done by multiplying by the **centering matrix**:

$$\mathbf{C}_{n \times n} = \mathbf{M}(1 - 1/n, -1/n) = \begin{bmatrix} 1 - 1/n & -1/n & \cdots & -1/n \\ -1/n & 1 - 1/n & \cdots & -1/n \\ \vdots & \vdots & \ddots & \vdots \\ -1/n & -1/n & \cdots & 1 - 1/n \end{bmatrix}.$$

This is a projection matrix with minor eigenvalue $\lambda_{**} = 0$ and major eigenvalue $\lambda_* = 1$. It is the annihilator corresponding to the projection onto the column space of the unit vector, so it has null space $\ker(\mathbf{C}) = \{\mathbf{a} \in \mathbb{R}^n | \mathbf{a} \propto \mathbf{1}\}$. When applied as a pre-multiplier to a column vector it transforms by subtracting the sample mean of the elements from all elements of the vector:

$$\mathbf{C}_{n\times n}\mathbf{x}_n = \begin{bmatrix} 1-1/n & -1/n & \cdots & -1/n \\ -1/n & 1-1/n & \cdots & -1/n \\ \vdots & \vdots & \ddots & \vdots \\ -1/n & -1/n & \cdots & 1-1/n \end{bmatrix} \begin{bmatrix} x_1 \\ x_2 \\ \vdots \\ x_n \end{bmatrix} = \begin{bmatrix} x_1 - \bar{x} \\ x_2 - \bar{x} \\ \vdots \\ x_n - \bar{x} \end{bmatrix}.$$

Centering can be extended to a design matrix composed of $m$ columns that are vectors to be centered around the row or column means. Pre-multiplication removes the column means (which is common in regression applications):

$$\mathbf{C}_{n\times n}\mathbf{x}_{n\times m} = \begin{bmatrix} x_{1,1} - \bar{x}_{\bullet 1} & x_{1,2} - \bar{x}_{\bullet 2} & \cdots & x_{1,m} - \bar{x}_{\bullet m} \\ x_{2,1} - \bar{x}_{\bullet 1} & x_{2,2} - \bar{x}_{\bullet 2} & \cdots & x_{2,m} - \bar{x}_{\bullet m} \\ \vdots & \vdots & \ddots & \vdots \\ x_{n,1} - \bar{x}_{\bullet 1} & x_{n,2} - \bar{x}_{\bullet 2} & \cdots & x_{n,m} - \bar{x}_{\bullet m} \end{bmatrix},$$

and post-multiplication removes the row means:

$$\mathbf{x}_{n\times m}\mathbf{C}_{m\times m} = \begin{bmatrix} x_{1,1} - \bar{x}_{1\bullet} & x_{1,2} - \bar{x}_{1\bullet} & \cdots & x_{1,m} - \bar{x}_{1\bullet} \\ x_{2,1} - \bar{x}_{2\bullet} & x_{2,2} - \bar{x}_{2\bullet} & \cdots & x_{2,m} - \bar{x}_{2\bullet} \\ \vdots & \vdots & \ddots & \vdots \\ x_{n,1} - \bar{x}_{n\bullet} & x_{n,2} - \bar{x}_{n\bullet} & \cdots & x_{n,m} - \bar{x}_{n\bullet} \end{bmatrix}.$$

Since the centering matrix can be used to remove row and column means, it can also be used to construct the class of all matrices with zero row and column sums. To do this, simple start with an arbitrary matrix $\mathbf{x}_{n\times m}$ and then remove both the row and column means to obtain:

$$\mathbf{C}_{n\times n}\mathbf{x}_{n\times m}\mathbf{C}_{m\times m} = \begin{bmatrix} x_{1,1} - \bar{x}_{1\bullet} - \bar{x}_{\bullet 1} + \bar{x}_{\bullet\bullet} & \cdots & x_{1,m} - \bar{x}_{1\bullet} - \bar{x}_{\bullet m} + \bar{x}_{\bullet\bullet} \\ x_{2,1} - \bar{x}_{2\bullet} - \bar{x}_{\bullet 1} + \bar{x}_{\bullet\bullet} & \cdots & x_{2,m} - \bar{x}_{2\bullet} - \bar{x}_{\bullet m} + \bar{x}_{\bullet\bullet} \\ \vdots & \ddots & \vdots \\ x_{n,1} - \bar{x}_{n\bullet} - \bar{x}_{\bullet 1} + \bar{x}_{\bullet\bullet} & \cdots & x_{n,m} - \bar{x}_{n\bullet} - \bar{x}_{\bullet m} + \bar{x}_{\bullet\bullet} \end{bmatrix}.$$

We have seen above that the centering matrix is the archetypal element of the class of singular non-degenerate double-constant matrices. The centering matrix is non-singular, so it cannot be inverted.[5] This is unsurprising, given that it is a projection which has the effect of removing the mean from a vector. Once the sample mean is lost, this information cannot be recovered endogenously —i.e., we require exogenous information to recover the lost mean. As with any other projection matrix, the centering matrix is idempotent, which is reflected in the fact that its eigenvalues are all equal to zero or one.

As with other double-constant matrices, the centering matrix gives a linear transform that can be viewed either as an operation in Euclidean space or in Fourier space. In the latter case we have already previously established that the Fourier values at the zero frequency are scaled by the minor eigenvalue and the Fourier values at all the non-zero frequencies are scaled by the major eigenvalue. In the case of the centering matrix this means that the Fourier values at the zero frequency are removed (scaled by an eigenvalue of zero) and the remaining Fourier values are left unchanged (scaled by an eigenvalue of one). This merely confirms a well-known aspect of Fourier analysis — removing the mean of a vector of values is equivalent to setting its power to zero at the zero frequency in Fourier space. The Plancherel theorem allows us to write the sums of squared deviations from the sample mean in Euclidean space or Fourier space as:

$$\sum_{i=1}^{n}(x_i - \bar{x})^2 = \|\boldsymbol{x} - \bar{\boldsymbol{x}}\|^2 = \boldsymbol{x}^{\mathrm{T}}\mathbf{C}_{n\times n}\boldsymbol{x} = \sum_{t=1}^{n-1}\|\mathcal{F}_{\boldsymbol{x}}(t)\|^2 \,.$$

---

[5] The Moore-Penrose pseudoinverse of the centering matrix is itself. This can easily be demonstrated directly by confirming that the matrix satisfies all criteria defining the pseudoinverse.

Parseval's theorem allows us to extend this to cross-products:

$$\sum_{i=1}^{n}(x_i - \bar{x})(y_i - \bar{y}) = (\boldsymbol{x} - \bar{\boldsymbol{x}}) \cdot (\boldsymbol{y} - \bar{\boldsymbol{y}}) = \boldsymbol{x}^{\mathrm{T}}\mathbf{C}_{n\times n}\boldsymbol{y} = \sum_{t=1}^{n-1}\mathcal{F}_{\boldsymbol{x}}(t) \cdot \overline{\mathcal{F}_{\boldsymbol{y}}(t)} \,.$$

These results are special cases of the general properties we previously derived for dot products of vectors transformed by arbitrary double-constant matrices. In these special cases the minor eigenvalue is zero and the major eigenvalue is one, so we obtain particularly simple forms composed of an unweighted sum of products of Fourier terms.

The centering matrix arises in a range of statistical applications dealing with random vectors. In such contexts it is common to use the centering matrix to simplify matrix expressions where the random vector is centered by removing its sample mean prior to some other operation. One such application occurs when calculating the Mahanalobis distance of a vector $\boldsymbol{x}$ using an $n \times n$ covariance matrix $\boldsymbol{\Sigma}$. The square of the Mahanalobis distance $D$ is defined as:

$$D^2(\boldsymbol{x}) = (\boldsymbol{x} - \bar{\boldsymbol{x}})^{\mathrm{T}}\boldsymbol{\Sigma}^{-1}(\boldsymbol{x} - \bar{\boldsymbol{x}}).$$

This squared-distance measure generalises the above sum-of-squares used in the Plancheral theorem. Using the spectral form $\boldsymbol{\Sigma} = \mathbf{A}\mathbf{V}\mathbf{A}^{\mathrm{T}}$ for the covariance matrix it can be written as:

$$D^2(\boldsymbol{x}) = (\boldsymbol{x} - \bar{\boldsymbol{x}})^{\mathrm{T}}\mathbf{A}\mathbf{V}^{-1}\mathbf{A}^{\mathrm{T}}(\boldsymbol{x} - \bar{\boldsymbol{x}}).$$

Using the centering matrix we can write the squared-distance as:

$$D(\boldsymbol{x}) = \boldsymbol{x}_{\mathrm{c}}\boldsymbol{\Sigma}^{-1}\boldsymbol{x}_{\mathrm{c}} = \boldsymbol{x}^{\mathrm{T}}\mathbf{A}_{\mathrm{c}}\mathbf{V}^{-1}\mathbf{A}_{\mathrm{c}}^{\mathrm{T}}\boldsymbol{x} \qquad \boldsymbol{x}_{\mathrm{c}} = \mathbf{C}_{n\times n}\boldsymbol{x} \quad \mathbf{A}_{\mathrm{c}} = \mathbf{C}_{n\times n}\mathbf{A}.$$

We see here that the squared Mahanalobis distance can be written as a quadratic form either for the centred data or on the original data using a centered eigenvalue matrix. This equivalent provides an extension to the Plancherel theorem above, but unlike in that theorem, the relevant eigenspace in the latter form is no longer Fourier space. (It is now a more complex eigenspace that is determined by the centered eigenvalue matrix of the covariance matrix $\boldsymbol{\Sigma}$.)

The centering matrix is a projection matrix in its own right, but it also arises in relation to more general projections occurring in linear models that include an intercept term. To see this, we consider the design matrix $x$ consisting of a column of ones (for an intercept term in the model) and a set of explanatory vectors $x_1, \ldots, x_m$. We denote the design matrix and its parts by:

$$x = [\mathbf{1}_n \quad x_{n \times m}] = [\mathbf{1}_n \quad x_1 \quad \cdots \quad x_m].$$

This matrix is used in the multiple linear regression model $y = x\beta + \varepsilon$, where projection onto the column space of the design matrix is done using the projection matrix $\mathbf{H}_{n \times n} = x(x^T x)^{-1} x^T$ (more commonly known as the "hat matrix"). In this context it is often useful to decouple the estimation of the intercept term from the estimation of the coefficients of the linear transform applying to the explanatory variables, so the annihilator matrix $I - \mathbf{H}_{n \times n}$ is decomposed into a first part centering the data and another part projecting onto the column space of the *centered* explanatory vectors. If we define $x_c = \mathbf{C}_{n \times n} x_{n \times m}$ as the centered explanatory data then the overall annihilator matrix can be decomposed via the blockwise formula for projections as:

$$I - \mathbf{H}_{n \times n} = \mathbf{C}_{n \times n} - x_c (x_c^T x_c)^{-1} x_c^T.$$

Applying this projection gives residual vector $r = (I - \mathbf{H}_{n \times n}) y = \mathbf{C}_{n \times n} y - x_c (x_c^T x_c)^{-1} x_c^T y$, which consists of the centered response vector (showing deviations of the response variables from their sample mean) plus predictive deviations that are formed by projecting the response vector onto the column space of the centered explanatory data. Preliminary centering of data in regression analysis may be undertaken for interpretive reasons, of for geometric analysis, though it has sometimes also been applied in the mistaken belief that it alleviates problems with collinearity in regressions using interaction terms (see Echambadi and Hess 2007).

Finally, we can gain some insight into ordinary least squares (OLS) estimation by looking at it in terms of the data. The estimated coefficient vector $\hat{\beta} = (x_c^T x_c)^{-1} x_c^T y$ from OLS estimation satisfies the normal equations $(x_c^T x_c) \hat{\beta} = x_c^T y$. In scalar form these linear equations are:

$$\sum_{k=1}^{m} [x_c^T x_c]_{i,k} \hat{\beta}_i = [x_c^T y]_i \qquad \text{for } i = 1, \ldots n.$$

The terms of both the bracketed matrices in this expression can be represented in Fourier space using Parseval's theorem. Taking $\mathcal{F}_x(t) = [\mathcal{F}_{x_1}(t) \quad \cdots \quad \mathcal{F}_{x_m}(t)]$ as the Fourier representation for all of the explanatory vectors, and using $\circ$ to denote the Hadamard (elementwise) product, we can write these normal equations in alternative form as:

$$\sum_{t=1}^{n-1} (\hat{\beta} \circ \mathcal{F}_x(t)) \cdot \overline{\mathcal{F}_x(t)} = \sum_{t=1}^{n-1} \mathcal{F}_x(t) \cdot \overline{\mathcal{F}_y(t)}.$$

## 4. The centering matrix and degrees-of-freedom

One additional use of the centering matrix in statistics is to assist in establishing the meaning of "degrees-of-freedom" and showing its connection to estimation of variance parameters in statistical models. In statistical models, it is common to estimate variance parameters with statistics akin to the sample variance, using Bessel's correction in the denominator to offset the effect of having estimated the mean from the sample data. In linear models, the denominator in the variance estimator is equal to the residual degrees-of-freedom in the model, but it is not immediately clear why this is the case.[6]

Here we illustrate a simple example where we have an IID sample $x = (x_1, \ldots, x_n)$ taken from a fixed distribution with mean $\mathbb{E}(X) = \mu$ and variance $\mathbb{V}(X) = \sigma^2 < \infty$. Suppose we estimate the mean using the sample mean, which gives us a vector of deviations from the mean, which is $r = x - \bar{x} = \mathbf{C}_{n \times n} x$. Formally, the **degrees-of-freedom** is the dimension of the space of possible values for the residual vector $r \in \mathcal{R}$, which is:

$$DF \equiv \dim \mathcal{R} = \dim\{\mathbf{C}_{n \times n} x | x \in \mathbb{R}^n\} = \operatorname{rank} \mathbf{C}_{n \times n}.$$

Now, for IID data we have $\mathbb{E}(\|\mathcal{F}_X(t)\|^2) = \sigma^2$ over all the discrete frequency values $t$. Hence, using the eigendecomposition of the centering matrix (and labelling its eigenvalues without substitution of their actual values for now), and recognising that $\operatorname{tr} \mathbf{C}_{n \times n} = \operatorname{rank} \mathbf{C}_{n \times n}$ since the centering matrix is a projection matrix, we have:

$$\mathbb{E}(\|X - \bar{X}\|^2) = \mathbb{E}(X^\mathrm{T} \mathbf{C}_{n \times n} X) = \mathbb{E}\left(\sum_{t=0}^{n-1} \lambda_i \cdot \|\mathcal{F}_X(t)\|^2\right)$$

$$= \sum_{t=0}^{n-1} \lambda_i \cdot \mathbb{E}(\|\mathcal{F}_X(t)\|^2)$$

$$= \sigma^2 \sum_{t=0}^{n-1} \lambda_i$$

$$= \sigma^2 \operatorname{tr} \mathbf{C}_{n \times n}$$

$$= \sigma^2 \operatorname{rank} \mathbf{C}_{n \times n} = \sigma^2 \times DF.$$

---

[6] Heuristically, this occurs because deviations from an estimated mean are always smaller than deviations from the true mean; division by the residual degrees-of-freedom, instead of the number of observations, inflates the estimator to account for this. This heuristic explanation is generally well-understood, but the *formal* connection from bias-correction to the degrees-of-freedom is not widely known.

The standard (unbiased) sample variance estimator is thus formed as:

$$S^2 = \frac{\|X - \bar{X}\|^2}{DF} = \frac{1}{n-1} \sum_{i=1}^{n} (X_i - \bar{X})^2.$$

Division by $DF$ (rather than by $n$) in this expression is known as "Bessel's correction". It is a result that famously taunts undergraduate statistics students, and often leaves their instructors flailing with heuristic explanations. Our analysis establishes a connection between Bessel's correction and the degrees-of-freedom, related through the eigenvalues of the centering matrix.

The above result is for a simple case where we have no explanatory variables, but it can easily be extended to unbiased estimation of the error variance in linear models where we have one or more explanatory variables. In all these cases, Bessel's correction and the degrees-of-freedom are connected through the eigenvalues of the centering matrix. The eigenvalues values directly determine the degrees-of-freedom through the dimension of the space of allowable centered vectors. The eigenvalues also directly affect the expected value of the squared-norm of the centered vectors, which is the basis for estimation of the variance parameter.

## 5. Decomposition of double-constant and centering matrices

In some applications it can be useful to decompose the double-constant matrix into blocks that correspond to some partition of a sample data vector. To do this, suppose we decompose the matrix by taking some dimension values $n_1 + n_2 + \cdots + n_k = n$. We can easily establish the general block decomposition:

$$\mathbf{M}(a, t) = \begin{bmatrix} \mathbf{M}_1(a,t) & t \cdot \mathbf{1}_{n_1 \times n_2} & \cdots & t \cdot \mathbf{1}_{n_1 \times n_{k-1}} & t \cdot \mathbf{1}_{n_1 \times n_k} \\ t \cdot \mathbf{1}_{n_2 \times n_1} & \mathbf{M}_2(a,t) & \cdots & t \cdot \mathbf{1}_{n_2 \times n_{k-1}} & t \cdot \mathbf{1}_{n_2 \times n_k} \\ \vdots & \vdots & \ddots & \vdots & \vdots \\ t \cdot \mathbf{1}_{n_{k-1} \times n_1} & t \cdot \mathbf{1}_{n_{k-1} \times n_2} & \cdots & \mathbf{M}_{k-1}(a,t) & t \cdot \mathbf{1}_{n_{k-1} \times n_k} \\ t \cdot \mathbf{1}_{n_k \times n_1} & t \cdot \mathbf{1}_{n_k \times n_2} & \cdots & t \cdot \mathbf{1}_{n_k \times n_{k-1}} & \mathbf{M}_k(a,t) \end{bmatrix}$$

$$= \begin{bmatrix} \mathbf{M}_1(a_1, t_1) & 0 & \cdots & 0 \\ 0 & \mathbf{M}_2(a_2, t_2) & \cdots & 0 \\ \vdots & \vdots & \ddots & \vdots \\ 0 & 0 & \cdots & \mathbf{M}_k(a_k, t_k) \end{bmatrix}$$

$$+ \begin{bmatrix} \mathbf{M}_1(a - a_1, t - t_1) & t \cdot \mathbf{1}_{n_1 \times n_2} & \cdots & t \cdot \mathbf{1}_{n_1 \times n_k} \\ t \cdot \mathbf{1}_{n_2 \times n_1} & \mathbf{M}_2(a - a_2, t - t_2) & \cdots & t \cdot \mathbf{1}_{n_2 \times n_k} \\ \vdots & \vdots & \ddots & \vdots \\ t \cdot \mathbf{1}_{n_k \times n_1} & t \cdot \mathbf{1}_{n_k \times n_2} & \cdots & \mathbf{M}_k(a - a_k, t - t_k) \end{bmatrix}.$$

This decomposition expresses the $n \times n$ double-constant matrix as a block matrix composed of smaller double-constant matrices on the diagonals, plus a remainder term that is also a block matrix composed of double-constant matrices and constant matrices. (Note that the remainder is not generally a double-constant matrix itself.)

One area where this form of decomposition is useful is when we are centering data vectors and we wish to break our analysis down to examine a partition of the data. This leads to situations where it is useful to decompose the centering matrix. In order to apply this decomposition to the centering matrix we use the values:

$$a = 1 - 1/n \qquad t = -1/n,$$
$$a_\ell = 1 - 1/n_\ell \qquad t_\ell = -1/n_\ell.$$

Taking these values gives:

$$\mathbf{M}_\ell(a - a_\ell, t - t_\ell) = \frac{1}{n} \cdot w_\ell \cdot \mathbf{1}_{n_\ell \times n_\ell} \qquad w_\ell \equiv \frac{n - n_\ell}{n_\ell}.$$

Hence, we obtain the following decomposition for the centering matrix:

$$\mathbf{C}_n = \begin{bmatrix} \mathbf{C}_1 & \mathbf{0} & \cdots & \mathbf{0} \\ \mathbf{0} & \mathbf{C}_2 & \cdots & \mathbf{0} \\ \vdots & \vdots & \ddots & \vdots \\ \mathbf{0} & \mathbf{0} & \cdots & \mathbf{C}_k \end{bmatrix} + \frac{1}{n} \begin{bmatrix} w_1 \cdot \mathbf{1}_{n_1 \times n_1} & -\mathbf{1}_{n_1 \times n_2} & \cdots & -\mathbf{1}_{n_1 \times n_k} \\ -\mathbf{1}_{n_2 \times n_1} & w_2 \cdot \mathbf{1}_{n_2 \times n_2} & \cdots & -\mathbf{1}_{n_2 \times n_k} \\ \vdots & \vdots & \ddots & \vdots \\ -\mathbf{1}_{n_k \times n_1} & -\mathbf{1}_{n_k \times n_2} & \cdots & w_k \cdot \mathbf{1}_{n_k \times n_k} \end{bmatrix}.$$

This decomposition is useful to obtain a corresponding decomposition for the sums-of-squares quantities for partitioned sample. Consider a pooled data vector $\mathbf{x} = (x_1, \ldots, x_n)$ partitioned into subgroups $\mathbf{x} = (\mathbf{x}_1, \mathbf{x}_2, \ldots, \mathbf{x}_k)$ with respective sample sizes $n_1, n_2, \ldots, n_k$. The sums-of-squares of the pooled data and the subgroups can be written using the centering matrix as:

$$\text{SS}_\text{p} = \mathbf{x}^\text{T} \mathbf{C} \mathbf{x} \qquad \text{SS}_\ell = \mathbf{x}_\ell^\text{T} \mathbf{C}_\ell \mathbf{x}_\ell.$$

Hence, we can obtain the sum-of-squares decomposition:

$$\text{SS}_\text{p} = \mathbf{x}^\text{T} \mathbf{C} \mathbf{x} = \sum_{\ell=1}^{k} \mathbf{x}_\ell^\text{T} \mathbf{C}_\ell \mathbf{x}_\ell + \frac{1}{n} \cdot \mathbf{x}^\text{T} \begin{bmatrix} w_1 \cdot \mathbf{1}_{n_1 \times n_1} & -\mathbf{1}_{n_1 \times n_2} & \cdots & -\mathbf{1}_{n_1 \times n_k} \\ -\mathbf{1}_{n_2 \times n_1} & w_2 \cdot \mathbf{1}_{n_2 \times n_2} & \cdots & -\mathbf{1}_{n_2 \times n_k} \\ \vdots & \vdots & \ddots & \vdots \\ -\mathbf{1}_{n_k \times n_1} & -\mathbf{1}_{n_k \times n_2} & \cdots & w_k \cdot \mathbf{1}_{n_k \times n_k} \end{bmatrix} \mathbf{x}$$

$$= \sum_{\ell=1}^{k} \text{SS}_\ell + \frac{1}{n} \left[ \sum_\ell w_\ell \cdot \mathbf{x}_\ell^\text{T} \mathbf{1}_{n_\ell \times n_\ell} \mathbf{x}_\ell - \sum_{\ell \neq h} \mathbf{x}_\ell^\text{T} \mathbf{1}_{n_\ell \times n_h} \mathbf{x}_h \right]$$

$$= \sum_{\ell=1}^{k} \text{SS}_\ell + \frac{1}{n} \left[ \sum_\ell w_\ell \cdot \sum_{i=1}^{n_\ell} \sum_{j=1}^{n_\ell} x_{\ell,i} x_{\ell,j} - \sum_{\ell \neq h} \sum_{i=1}^{n_\ell} \sum_{j=1}^{n_\ell} x_{\ell,i} x_{h,j} \right].$$

In the case where $k = 2$ (i.e., a pooled sample composed of two subgroups) the weightings reduce to $w_\ell = n_1 n_2 / n_\ell^2$ and so we get the simpler decomposition:

$$\begin{aligned}
SS_p &= SS_1 + SS_2 + \frac{1}{n} \cdot \mathbf{x}^T \begin{bmatrix} w_1 \cdot \mathbf{1}_{n_1 \times n_1} & -\mathbf{1}_{n_1 \times n_2} \\ -\mathbf{1}_{n_2 \times n_1} & w_2 \cdot \mathbf{1}_{n_2 \times n_2} \end{bmatrix} \mathbf{x} \\
&= SS_1 + SS_2 + \frac{n_1 n_2}{n_1 + n_2} \cdot [\mathbf{x}_1^T \ \mathbf{x}_2^T] \begin{bmatrix} \mathbf{1}_{n_1 \times n_1}/n_1^2 & -\mathbf{1}_{n_2 \times n_1}/n_1 n_2 \\ -\mathbf{1}_{n_2 \times n_1}/n_1 n_2 & \mathbf{1}_{n_1 \times n_1}/n_2^2 \end{bmatrix} \begin{bmatrix} \mathbf{x}_1 \\ \mathbf{x}_2 \end{bmatrix} \\
&= SS_1 + SS_2 + \frac{n_1 n_2}{n_1 + n_2} \cdot [\mathbf{x}_1^T \ \mathbf{x}_2^T] \begin{bmatrix} \mathbf{1}_{n_1}/n_1 \\ -\mathbf{1}_{n_2}/n_2 \end{bmatrix} \begin{bmatrix} \mathbf{1}_{n_1}/n_1 \\ -\mathbf{1}_{n_2}/n_2 \end{bmatrix}^T \begin{bmatrix} \mathbf{x}_1 \\ \mathbf{x}_2 \end{bmatrix} \\
&= SS_1 + SS_2 + \frac{n_1 n_2}{n_1 + n_2} \cdot (\bar{x}_1 - \bar{x}_2)^2.
\end{aligned}$$

This is a useful decomposition result that can be used when dealing with situations where a pooled sample is obtained from two smaller subgroups. (This sum-of-squares decomposition leads easily to a corresponding decomposition for the sample variance given in O'Neill 2014, Result 1, pp. 283.)

### 6. The equicorrelation matrix and its use in statistical models

Another important class of double-constant matrices are the positive definite forms. Every positive-definite double-constant matrix is a positive multiple of the **equicorrelation matrix**, which is an $n \times n$ double-constant matrix given by:[7]

$$\mathbf{\Sigma}_\rho^2 = \mathbf{M}(1, \rho) = \begin{bmatrix} 1 & \rho & \cdots & \rho \\ \rho & 1 & \cdots & \rho \\ \vdots & \vdots & \ddots & \vdots \\ \rho & \rho & \cdots & 1 \end{bmatrix} \qquad -\frac{1}{n-1} < \rho < 1.$$

Since the equicorrelation matrix is a correlation matrix, it is positive definite, which means the correlation parameter must fall within the specified range (equivalent to requiring that both the major and minor eigenvalues are positive). When $\rho \to 1$ the equicorrelation matrix approaches the unit matrix and when $\rho \to -1/(n-1)$ it approaches a multiple of the centering matrix. If we wish to allow extension of this matrix for arbitrary $n$ then we can impose the stronger requirement that $0 \leq \rho < 1$, which ensures that the minor eigenvalues is positive for all $n$.

---

[7] The reader may have noticed that our notation for the equicorrelation matrix uses a square superscript. The reason for this is that the matrix represents a variance matrix, and we will want to use its principal square-root, the **equicorrelation root-matrix** $\mathbf{\Sigma}_\rho$ in some of our later discussion.

The major and minor eigenvalues of the equicorrelation matrix are given by $\lambda_* = 1 - \rho$ and $\lambda_{**} = 1 - \rho + n\rho$ respectively. The equicorrelation matrix and related matrices can be written in canonical form (i.e., in terms of the major and minor eigenvalues) as:

$$\mathbf{\Sigma}_\rho^2 = \lambda_* \cdot \mathbf{I}_{n \times n} + \frac{\lambda_{**} - \lambda_*}{n} \cdot \mathbf{1}_{n \times n} \qquad \mathbf{\Sigma}_\rho = \sqrt{\lambda_*} \cdot \mathbf{I}_{n \times n} + \frac{\sqrt{\lambda_{**}} - \sqrt{\lambda_*}}{n} \cdot \mathbf{1}_{n \times n},$$

$$\mathbf{\Sigma}_\rho^{-2} = \frac{1}{\lambda_*} \cdot \mathbf{I}_{n \times n} + \frac{1}{n}\left(\frac{1}{\lambda_{**}} - \frac{1}{\lambda_*}\right) \cdot \mathbf{1}_{n \times n} \qquad \mathbf{\Sigma}_\rho^{-1} = \frac{1}{\sqrt{\lambda_*}} \cdot \mathbf{I}_{n \times n} + \frac{1}{n}\left(\frac{1}{\sqrt{\lambda_{**}}} - \frac{1}{\sqrt{\lambda_*}}\right) \cdot \mathbf{1}_{n \times n}.$$

We can rewrite these in explicit form as:

$$\mathbf{\Sigma}_\rho^2 = (1 - \rho) \cdot \mathbf{I}_{n \times n} + \rho \cdot \mathbf{1}_{n \times n},$$

$$\mathbf{\Sigma}_\rho^{-2} = \frac{1}{1 - \rho} \cdot \mathbf{I}_{n \times n} - \frac{\rho}{(1 - \rho)(1 - \rho + n\rho)} \cdot \mathbf{1}_{n \times n},$$

$$\mathbf{\Sigma}_\rho = \sqrt{1 - \rho} \cdot \mathbf{I}_{n \times n} + \frac{\sqrt{1 - \rho + n\rho} - \sqrt{1 - \rho}}{n} \cdot \mathbf{1}_{n \times n},$$

$$\mathbf{\Sigma}_\rho^{-1} = \frac{1}{\sqrt{1 - \rho}} \cdot \mathbf{I}_{n \times n} + \frac{1}{n}\left[\frac{1}{\sqrt{1 - \rho + n\rho}} - \frac{1}{\sqrt{1 - \rho}}\right] \cdot \mathbf{1}_{n \times n}.$$

In statistical models with equicorrelated errors the variance matrix of the error vector will have off-diagonal elements that are restricted in such a way as to yield equal correlation between each pair of error terms. Construction of variance matrices of this form will make use of the equicorrelation matrix. Early statistical literature on models of this form looked at the effects of equicorrelation on statistical tests (Walsh 1946, Halperin 1951, Geisser 1963, Srivastava 1965) and later literature constructed tests of the correlation coefficient of equicorrelated data (Sengupta 1987, Viana 1994, Bhatti and Barry 1996). Statistical models involving equally correlated error terms have been applied regularly in economics and finance (Elton and Gruber 1973, Berndt and Savin 1975, Engle and Kelly 2012, Aboura and Chevallier 2014, Clements *et al* 2015, O'Neill 2015, Kurose and Omori 2016).

In statistical problems involving estimation of an unknown variance parameter, equicorrelation between the observable values is sometimes regarded as changing the "effective sample size", where the latter is regarded as being the sample size that would be equivalent in information if the values were uncorrelated. Thus, if the observable values are positively correlated, then they give less information about the variance than if they were uncorrelated, and so the "effective" sample size is lower than the sample size. Conversely, if the observable values are negatively correlated then they give more information about the variance than if they were uncorrelated, and so the "effective" sample size is higher than the sample size.

To illustrate this phenomenon, we will consider a simple example where we have an observed sample $x = (x_1, ..., x_n)$ from a fixed marginal distribution (with finite variance), but with an allowance for non-zero equicorrelation between the elements. This implies that we have an observable vector $X$ with mean vector and variance matrix given by:
$$\mathbb{E}(X) = \mu \mathbf{1}_n \qquad \mathbb{V}(X) = \sigma^2 \mathbf{\Sigma}_\rho^2.$$
We can decompose the observable random vector as $X = \mu \mathbf{1}_n + \sigma \mathbf{\Sigma}_\rho \epsilon$ where we have a vector of underlying uncorrelated error terms with mean vector and variance matrix given by:
$$\mathbb{E}(\epsilon) = \mathbf{0} \qquad \mathbb{V}(\epsilon) = \mathbf{I}_{n \times n}.$$
Applying the multiplication rule in Theorem 5, it can be shown that $\mathbf{C}_{n \times n} \mathbf{\Sigma}_\rho = (1 - \rho) \mathbf{C}_{n \times n}$. Thus, the residual deviation vector for the sample can be written as:
$$\begin{aligned} X - \bar{X} = \mathbf{C}_{n \times n} X &= \mathbf{C}_{n \times n} (\mu \mathbf{1}_n + \sigma \mathbf{\Sigma}_\rho \epsilon) \\ &= \mathbf{0}_{n \times n} + \sigma \mathbf{C}_{n \times n} \mathbf{\Sigma}_\rho \epsilon \\ &= \sigma \mathbf{C}_{n \times n} \mathbf{\Sigma}_\rho \epsilon \\ &= \sigma (1 - \rho) \mathbf{C}_{n \times n} \epsilon. \end{aligned}$$
From this expression we see that the variance and equicorrelation parameters affect the residual deviation vector only through the product $\sigma(1 - \rho)$.

If the parameter $\rho$ is known then we can use this to adjust the estimator for the variance. In this case, we will show that the "effective sample size" is $n_{\text{Eff}} = 1 + (n-1)(1-\rho)^2$, so we have corresponding "effective" degrees-of-freedom $DF_{\text{Eff}} = n_{\text{Eff}} - 1$. The latter value relates to the standard degrees-of-freedom by the equation:
$$DF_{\text{Eff}} = (1 - \rho)^2 DF.$$
We will now assume that the error terms are IID, which is a stronger assumption than that they are merely uncorrelated. Under this assumption, we have $\mathbb{E}(\|\mathcal{F}_\epsilon(t)\|^2) = 1$ over all the discrete frequency values $t$. Applying the above form of the residual deviation vector we have:
$$\begin{aligned} \mathbb{E}(\|X - \bar{X}\|^2) &= \sigma^2 (1 - \rho)^2 \mathbb{E}(\epsilon^{\mathsf{T}} \mathbf{C}_{n \times n} \epsilon) \\ &= \sigma^2 (1 - \rho)^2 \operatorname{rank} \mathbf{C}_{n \times n} \\ &= \sigma^2 (1 - \rho)^2 \times DF \\ &= \sigma^2 \times DF_{\text{Eff}}. \end{aligned}$$

Thus, if $\rho$ is known, then the standard (unbiased) sample variance estimator is formed as:

$$S^2 = \frac{\|X - \bar{X}\|^2}{DF_{\text{Eff}}} = \frac{1}{n_{\text{Eff}} - 1} \sum_{i=1}^{n} (X_i - \bar{X})^2.$$

This confirms the earlier assertion that equicorrelation changes the "effective" sample size in problems involving estimation of variance. If $\rho = 0$ the observable values are uncorrelated and we have $n_{\text{Eff}} = n$, which gives the standard sample variance estimator. If $\rho > 0$ we have $1 \leq n_{\text{Eff}} < n$ and if $\rho < 0$ we have $n < n_{\text{Eff}} < n + (2n-1)/(n-1)$.[8]

The above mathematics shows that the equicorrelation parameter affects the "effective" sample size for estimation of the variance. The estimation results we have given assumes a context where the equicorrelation parameter is known, which is not very realistic. In practice, the equicorrelation parameter is usually unknown and it is necessary to estimate this parameter and the variance parameter separately (see e.g., Olkin and Pratt 1958, Donner and Bull 1983, Paul 1990, Paul and Barnwal 1990, Viana 1994). The latter problem requires us to observe multiple sample vectors $x_1, \ldots, x_N$, so we use multiple sample vectors to estimate the mean parameter, and this allows us to disentangle the variance parameter from the equicorrelation parameter. Broader problems and models involving equicorrelation are beyond the scope of our analysis, but we note that the equicorrelation matrix appears in all these broader problems, and it often interacts with the centering matrix in problems that also involve mean estimation.

## 7. Indefinite and negative-definite forms

The non-degenerate form of the double constant matrix occurs when we have $a \neq t \neq 0$. For this class of matrices, we have seen that the singular form of the matrix is proportionate to the centering matrix and the positive-definite form is proportionate to the equicorrelation matrix. The remaining cases are the indefinite and negative-definite forms, which we examine here. These last forms occur in the following cases:

| | | | | |
|---|---|---|---|---|
| $t < 0$ | $0 < a - t < -nt$ | $\lambda_* > 0$ | $\lambda_{**} < 0$ | Indefinite, |
| $t > 0$ | $-nt < a - t < 0$ | $\lambda_* < 0$ | $\lambda_{**} > 0$ | Indefinite, |
| $t < 0$ | $a - t < 0$ | $\lambda_* < 0$ | $\lambda_{**} < 0$ | Negative Definite, |
| $t > 0$ | $a - t < -nt$ | $\lambda_* < 0$ | $\lambda_{**} < 0$ | Negative Definite. |

---

[8] We assume here that $n > 1$ so that there is a basis for variance estimation.

In the first section to this paper we defined the double-constant matrix as a weighted sum of the identity matrix and the unit matrix. We also gave a canonical form for the double-constant matrix, by writing this decomposition with weights that were the eigenvalues. The first form is a decomposition of the identity matrix and the unit matrix, which is a natural decomposition, insofar as it splits the double-constant matrix into a weighted sum of two simple matrices with well-known properties. The canonical form is also very useful, and we will now have a more detailed look at this form. In fact, both of these decompositions are just special cases of the more general rule that any double-constant matrix can be written as a weighted sum of *any* two double-constant matrices that are not proportional to each other. We give this general rule in Theorem 7 and we write the canonical form as Theorem 8.

**THEOREM 7:** Let $\mathbf{M}_1$ and $\mathbf{M}_2$ be any two double-constant matrices that are not proportional to each other. Then for any $a, t \in \mathbb{R}$ there exist values $A(a,t)$ and $T(a,t)$ such that:
$$\mathbf{M}(a,t) = A(a,t) \cdot \mathbf{M}_1 + T(a,t) \cdot \mathbf{M}_2.$$

**THEOREM 8:** For any $a, t \in \mathbb{R}$ we have the **canonical form**:
$$\mathbf{M}(a,t) = \lambda_* \cdot \mathbf{C}_{n \times n} + \lambda_{**} \cdot \frac{\mathbf{1}_{n \times n}}{n},$$
where $\lambda_* = a - t$ and $\lambda_{**} = a - t + nt$ are the major and minor eigenvalues respectively.

**COROLLARY:** Letting $\mathbf{\Sigma}_*^2$ and $\mathbf{\Sigma}_{**}^2$ denote the respective lower and upper limits (i.e., taking the limits $\rho \to -1/(n-1)$ and $\rho \to 1$ respectively) of the equicorrelation matrix, we have:
$$\mathbf{M}(a,t) = \frac{n-1}{n} \cdot \lambda_* \cdot \mathbf{\Sigma}_*^2 + \frac{1}{n} \cdot \lambda_{**} \cdot \mathbf{\Sigma}_{**}^2.$$

The decomposition in Theorem 8 is a useful alternative to the initial decomposition we used as the defining form of the double-constant matrix. The decomposition in the canonical form has weightings given by the eigenvalues of the matrix, so it is particularly useful for characterising double-constant matrices in terms of their "definiteness" conditions. We can see that each of the four sets of conditions above —for the indefinite and negative-definite forms— correspond to combinations of positive and negative weights on the centering matrix $\mathbf{C}_{n \times n}$ and the mean matrix $\mathbf{1}_{n \times n}/n$. The corollary frames the decomposition slightly differently, in terms of the limits of the equicorrelation matrix. In this latter form the double-constant matrix is a convex combination of the eigenvalue-weighted limits of the limiting equicorrelation matrices, with the convex combination being proportional to the multiplicity of the eigenvalues.

Theorem 8 and its corollary establish the centering matrix and equicorrelation matrix as being especially important for describing the class of double-constant matrices. We can write every double-constant matrix as a weighted sum of lower and upper limiting equicorrelation matrices, with the signs of the weights corresponding to the signs of the major and minor eigenvalues. Alternatively, we can write every double-constant matrix as a weighted sum of lower and upper limiting equicorrelation matrices, with the weights corresponding directly to the major and minor eigenvalues (without any further weighting). Thus, indefinite forms occur when one of these weights is positive and the other is negative, and negative-definite forms occur when both weights are negative. This result brings us full-circle, back to the initial characterisation of the class of double-constant matrices. By examining the centering matrix and equicorrelation matrix as special cases, we are able to concentrate on those matrices that arise most commonly in statistical applications, and also give a canonical characterisation of the matrix class that is most useful for determining the definiteness of the matrix.

## 8. Summary and final remarks

In this paper we have examined double-constant matrices, which are square matrices composed of a constant main diagonal, and another constant for all off-diagonal elements, where these two constants may differ. This matrix is a special case of the circulant matrix, which is itself a special case of the Toeplitz matrix. Although constituting a special case of a wider matrix class, it is an interesting matrix form in its own right, manifesting in several important matrices that arise in statistical applications. We have seen that it is most useful to characterise double-constant matrices according to their major and minor eigenvalues, which yields a canonical form where the matrix is written as a weighted sum of the centering matrix and the mean matrix. In Figure 1 we partitioned the class of all double-constant matrices into important subgroups. Among the non-degenerate double-constant matrices (i.e., those that do not degenerate either to a constant matrix or a scalar multiple of the identity matrix) we identified that the singular form is proportional to the centering matrix and the positive-definite form is (positively) proportional to the equicorrelation matrix.

Our analysis of the double-constant matrix has looked at its general properties, including the eigendecomposition of the matrix, its relationship to Fourier analysis, and its closure under standard matrix operations. We have also looked at the specific properties of the centering matrix and equicorrelation matrix, and we have given a basic outline of how these matrices arise in statistical applications. We have used these matrices to elucidate certain statistical phenomena, such as the connection between degrees-of-freedom and Bessel's correction in variance estimation. By characterising the class of double-constant matrices in a canonical form that is closely connected to specific matrices used in statistical applications, we hope we have whet the reader's appetite for exploring the interesting ways in which statistical analysis closely relates to the spectral properties of matrices.

# Appendix – Proof of Theorems

In this appendix we will prove the theorems shown in the body of the paper plus and associated lemma relating to the discrete Fourier transformation.

**PROOF OF THEOREM 1:** For the case where $a \neq t$ we can apply the matrix determinant lemma to obtain:

$$\begin{aligned}
\det(\mathbf{M}) &= \det((a-t) \cdot \mathbf{I}_{n \times n} + t \cdot \mathbf{1}_{n \times n}) \\
&= \det((a-t) \cdot \mathbf{I}_{n \times n} + t \cdot \mathbf{1}_n \mathbf{1}_n^\mathrm{T}) \\
&= (1 + t \cdot \mathbf{1}_n^\mathrm{T}((a-t) \cdot \mathbf{I}_{n \times n})^{-1} \mathbf{1}_n) \det((a-t) \cdot \mathbf{I}) \\
&= \left(1 + \frac{t}{a-t} \cdot \mathbf{1}_n^\mathrm{T} \mathbf{1}_n\right) \det((a-t) \cdot \mathbf{I}_{n \times n}) \\
&= \left(1 + \frac{nt}{a-t}\right)(a-t)^n = (a-t)^{n-1}(a+t+nt).
\end{aligned}$$

The only remaining special case is $a = t$. Since $\det(\mathbf{1}_{n \times n}) = 0$ we have $\det(t \cdot \mathbf{I}_{n \times n}) = 0$ so the formula also holds in this case. ∎

**PROOF OF THEOREM 2:** The characteristic polynomial in this result can be obtained as a simple corollary to Theorem 1 by substituting $a - \lambda$ for $a$. This characteristic polynomial is already in factorised form, so this immediately establishes the eigenvalues. It remains only to find the corresponding eigenvectors. When $t = 0$ the double-constant matrix is proportional to the identity matrix (or the zero matrix if $a = 0$), so every non-zero vector is an eigenvector of each of its eigenvalues. In the non-trivial case when $t \neq 0$ we have:

$$\begin{aligned}
\mathbf{M} - \lambda_k \cdot \mathbf{I}_{n \times n} &= (a - t - \lambda_k) \cdot \mathbf{I}_{n \times n} + t \cdot \mathbf{1}_{n \times n} \\
&= \begin{cases} -nt \cdot \mathbf{I}_{n \times n} + t \cdot \mathbf{1}_{n \times n} & \text{for } k = 1, \\ t \cdot \mathbf{I}_{n \times n} & \text{for } k = 2, \ldots, n. \end{cases}
\end{aligned}$$

To facilitate our analysis we define $\dot{u}_k = \sum_{i=1}^n u_{k,i}$ as the sum of elements of the vector $\mathbf{u}_k$. The first eigenvector is found by setting:

$$\mathbf{0} = (\mathbf{M} - \lambda_1 \cdot \mathbf{I}_{n \times n})\mathbf{u}_1 = -(nt \cdot \mathbf{I}_{n \times n} - t \cdot \mathbf{1}_{n \times n})\mathbf{u}_1 = -t(n\mathbf{u}_1 - \dot{u}_1 \mathbf{1}_n),$$

so $\mathbf{u}_1 \propto \mathbf{1}_n$, which was to be shown. The remaining $n - 1$ eigenvectors are found by setting:

$$\mathbf{0} = (\mathbf{M} - \lambda_k \cdot \mathbf{I}_{n \times n})\mathbf{u}_k = t \cdot \mathbf{1}_{n \times n} \mathbf{u}_k = t\dot{u}_k \cdot \mathbf{1}_n,$$

so $\dot{u}_n = 0$, so any vector $\mathbf{u}_k$ with zero sum is an eigenvector corresponding to $\lambda_k$, which was to be shown. ∎

**LEMMA 1:** For all $r \in \mathbb{Z}$ we have:

$$\frac{1}{n}\sum_{k=0}^{n-1} \exp\left(-\frac{2\pi i r k}{n}\right) = \mathbb{I}(r \bmod n = 0).$$

**PROOF OF LEMMA 1:** Since the complex exponential is cyclic, we can restrict our attention to the integers $r = 0, 1, 2, \ldots, n-1$. For $r = 0$ we have:

$$\frac{1}{n}\sum_{k=0}^{n-1} \exp\left(-\frac{2\pi i r k}{n}\right) = \frac{1}{n}\sum_{k=0}^{n-1} \exp(0) = \frac{1}{n}\sum_{k=0}^{n-1} 1 = 1.$$

For the remaining values $r = 1, 2, \ldots, n-1$ we have $r/n \notin \mathbb{Z}$. Using the formula for the sum of a geometric progression we have:

$$\begin{aligned}
\frac{1}{n}\sum_{k=0}^{n-1} \exp\left(-\frac{2\pi i r k}{n}\right) &= \frac{1}{n}\sum_{k=0}^{n-1} \exp\left(-\frac{2\pi i r}{n}\right)^k \\
&= \frac{1}{n} \cdot \frac{1 - \exp(-2\pi i r)}{1 - \exp(-2\pi i r/n)} \\
&= \frac{1}{n} \cdot \frac{\exp(\pi i r/n)}{\exp(\pi i r/n)} \cdot \frac{1 - \exp(-2\pi i r)}{1 - \exp(-2\pi i r/n)} \\
&= \frac{1}{n} \cdot \frac{\exp(\pi i r/n)}{\exp(\pi i r)} \cdot \frac{\exp(\pi i r) - \exp(-\pi i r)}{\exp(\pi i r/n) - \exp(-\pi i r/n)} \\
&= \frac{1}{n} \cdot \frac{\exp(\pi i r/n)}{\exp(\pi i r)} \cdot \frac{\sin(\pi r)}{\sin(\pi r/n)} \\
&= \frac{1}{n} \cdot \frac{\exp(\pi i r/n)}{\exp(\pi i r)} \cdot \frac{0}{\sin(\pi r/n)} = 0.
\end{aligned}$$

This establishes the stated form for $r = 0, 1, 2, \ldots, n-1$ and the remaining values follow from the cyclic behaviour of the complex exponential. ∎

**PROOF OF THEOREM 3:** This theorem is actually a special case of Theorem 3.1 of Gray (2006), which establishes the eigendecomposition for circulant matrices. Notwithstanding this, we will prove the result here with respect only to double-constant matrices. From Theorem 2 we know that the columns $\mathbf{u}_1, \ldots, \mathbf{u}_n$ of $\mathbf{U}$ are eigenvectors for the double-constant matrix. The DFT $\mathbf{U}$ is a symmetric unitary square matrix, which means that $\mathbf{U}\bar{\mathbf{U}} = \bar{\mathbf{U}}\mathbf{U} = \mathbf{I}_{n \times n}$. We establish the diagonalization as follows. From Lemma 1 we have:

$$\frac{\mathbf{u}_r \cdot \mathbf{1}_n}{\sqrt{n}} = \frac{1}{n}\sum_{k=0}^{n-1} \exp\left(-\frac{2\pi i r k}{n}\right) = \mathbb{I}(r \bmod n = 0).$$

Similarly, for the corresponding rows of $\bar{\mathbf{U}}$ we also have:

$$\frac{\bar{\mathbf{u}}_r \cdot \mathbf{1}_n}{\sqrt{n}} = \frac{1}{n}\sum_{k=0}^{n-1} \exp\left(\frac{2\pi i r k}{n}\right) = \mathbb{I}(r \bmod n = 0).$$

Pre-multiplying $\mathbf{M}$ by $\bar{\mathbf{U}}$ therefore gives:

$$\begin{aligned}
\bar{\mathbf{U}}\mathbf{M} &= \bar{\mathbf{U}}\big((a-t)\cdot \mathbf{I}_{n\times n} + t\cdot \mathbf{1}_{n\times n}\big) \\
&= (a-t)\bar{\mathbf{U}}\mathbf{I}_{n\times n} + t\bar{\mathbf{U}}\mathbf{1}_{n\times n} \\
&= (a-t)\bar{\mathbf{U}}\mathbf{I}_{n\times n} + t\begin{bmatrix}\bar{\mathbf{u}}_1 \\ \bar{\mathbf{u}}_2 \\ \vdots \\ \bar{\mathbf{u}}_n\end{bmatrix}\begin{bmatrix}\mathbf{1}_n & \mathbf{1}_n & \cdots & \mathbf{1}_n\end{bmatrix} \\
&= (a-t)\bar{\mathbf{U}} + t\sqrt{n}\begin{bmatrix}\mathbf{1}_n^{\mathrm{T}} \\ \mathbf{0}_n^{\mathrm{T}} \\ \vdots \\ \mathbf{0}_n^{\mathrm{T}}\end{bmatrix}.
\end{aligned}$$

Post-multiplying this by $\mathbf{U}$ then gives:

$$\begin{aligned}
\bar{\mathbf{U}}\mathbf{M}\mathbf{U} &= (a-t)\bar{\mathbf{U}}\mathbf{U} + t\sqrt{n}\begin{bmatrix}\mathbf{1}_n^{\mathrm{T}} \\ \mathbf{0}_n^{\mathrm{T}} \\ \vdots \\ \mathbf{0}_n^{\mathrm{T}}\end{bmatrix}\begin{bmatrix}\mathbf{u}_1 & \cdots & \mathbf{u}_n\end{bmatrix} \\
&= (a-t)\mathbf{I}_{n\times n} + tn\begin{bmatrix}1 & 0 & \cdots & 0 \\ 0 & 0 & \cdots & 0 \\ \vdots & \vdots & \ddots & \vdots \\ 0 & 0 & \cdots & 0\end{bmatrix} \\
&= \begin{bmatrix}a-t+nt & 0 & \cdots & 0 \\ 0 & a-t & \cdots & 0 \\ \vdots & \vdots & \ddots & \vdots \\ 0 & 0 & \cdots & a-t\end{bmatrix} \\
&= \begin{bmatrix}\lambda_{**} & 0 & \cdots & 0 \\ 0 & \lambda_* & \cdots & 0 \\ \vdots & \vdots & \ddots & \vdots \\ 0 & 0 & \cdots & \lambda_*\end{bmatrix} = \boldsymbol{\Lambda}.
\end{aligned}$$

This establishes the first asserted diagonalization, and since $\mathbf{U}$ is orthonormal we also then have $\mathbf{U}\boldsymbol{\Lambda}\bar{\mathbf{U}} = \mathbf{U}(\bar{\mathbf{U}}\mathbf{M}\mathbf{U})\bar{\mathbf{U}} = (\mathbf{U}\bar{\mathbf{U}})\mathbf{M}(\mathbf{U}\bar{\mathbf{U}}) = \mathbf{M}$, which gives the corresponding eigendecomposition. The other diagonalisation $\mathbf{M} = \bar{\mathbf{U}}\boldsymbol{\Lambda}\mathbf{U}$ and the corresponding eigendecomposition follows from an analogous proof that is omitted here. It is worth noting here that the eigenvalues are given by the discrete Fourier transform of the vector of constant values $\mathbf{m} = (a, t, \ldots, t)$ taken from the double-constant matrix (removing the scaling constant in the unitary DFT). The minor eigenvalue is the Fourier transform at zero frequency and the major eigenvalue is the Fourier transform at each other frequency. ∎

**PROOF OF THEOREM 4:** Using the expanded form of the double-constant matrix we have:

$$\sum_{i=1}^{m} \kappa_i \cdot \mathbf{M}(a_i, t_i) = \sum_{i=1}^{m} \kappa_i[(a-t) \cdot \mathbf{I}_{n \times n} + t \cdot \mathbf{1}_{n \times n}]$$

$$= \sum_{i=1}^{m} \kappa_i(a-t) \cdot \mathbf{I}_{n \times n} + \sum_{i=1}^{m} \kappa_i t \cdot \mathbf{1}_{n \times n}$$

$$= \left(\sum_{i=1}^{m} \kappa_i(a-t)\right) \cdot \mathbf{I}_{n \times n} + \left(\sum_{i=1}^{m} \kappa_i t\right) \cdot \mathbf{1}_{n \times n}$$

$$= \mathbf{M}\left(\sum_{i=1}^{m} \kappa_i \cdot a_i, \sum_{i=1}^{m} \kappa_i \cdot t_i\right),$$

which was to be shown. ∎

**PROOF OF THEOREM 5:** The constants in the theorem correspond to the eigenvalues:

$$\lambda_{\times**} = \prod_{i=1}^{m}(a_i - t_i + nt_i) \qquad \lambda_{\times*} = \prod_{i=1}^{m}(a_i - t_i).$$

These are eigenvalues the products of the corresponding eigenvalues of the individual double-constant matrices in the product. Hence, using the eigendecomposition we have:

$$\prod_{i=1}^{m} \mathbf{M}(a_i, t_i) = \prod_{i=1}^{m} \mathbf{U}\mathbf{\Lambda}(a_i, t_i)\mathbf{\bar{U}} = \mathbf{U}\left[\prod_{i=1}^{m} \mathbf{\Lambda}(a_i, t_i)\right]\mathbf{\bar{U}} = \mathbf{U}\mathbf{\Lambda}_\times\mathbf{\bar{U}} = \mathbf{M}(a_\times, t_\times),$$

which was to be shown. ∎

**PROOF OF THEOREM 6:** To facilitate this proof we set $\lambda_k^* \equiv f(\lambda_k)$ to be updated eigenvalues and we use Sylvester's formula (REF) to obtain:

$$f(\mathbf{M}(a,t)) = \mathbf{U}f(\mathbf{\Lambda})\mathbf{\bar{U}}$$

$$= \mathbf{U}\,\text{diag}(f(\lambda_1), f(\lambda_2), \ldots, f(\lambda_n))\,\mathbf{\bar{U}}$$

$$= \mathbf{U}\,\text{diag}(\lambda_1^*, \lambda_2^*, \ldots, \lambda_n^*)\,\mathbf{\bar{U}}$$

$$= \mathbf{U}\mathbf{\Lambda}_*\mathbf{\bar{U}} = \mathbf{M}(a_*, t_*),$$

where:

$$a_* = \frac{1}{n}\sum_{k=1}^{n} \lambda_k^* = \frac{f(a-t+nt) + (n-1)f(a-t)}{n},$$

$$t_* = \frac{1}{n}(\lambda_1^* - \lambda_2^*) = \frac{f(a-t+nt) - f(a-t)}{n}.$$

The explicit form follows trivially from these values. ∎

**PROOF OF THEOREM 7:** Without loss of generality, let:

$$\mathbf{M}_1 = \mathbf{M}(a_1, t_1) \qquad \mathbf{M}_2 = \mathbf{M}(a_2, t_2).$$

Since $\mathbf{M}_1$ and $\mathbf{M}_2$ are not proportional to each other, we have $a_1 t_2 \neq a_2 t_1$. This allows us to unambiguously define the following functions:

$$A \equiv A(a,t) \equiv \frac{a t_2 - a_2 t}{a_1 t_2 - a_2 t_1} \qquad T \equiv T(a,t) \equiv \frac{a_1 t - a t_1}{a_1 t_2 - a_2 t_1},$$

which are the solutions to the simultaneous equations:

$$a = A a_1 + T a_2 \qquad t = A t_1 + T t_2.$$

Using these weights, we have:

$$\begin{aligned}
A \cdot \mathbf{M}_1 + T \cdot \mathbf{M}_2 &= A \cdot \mathbf{M}(a_1, t_1) + T \cdot \mathbf{M}(a_2, t_2) \\
&= \mathbf{M}(A a_1, A t_1) + \mathbf{M}(T a_2, T t_2) \\
&= \mathbf{M}(A a_1 + T a_2, A t_1 + T t_2) = \mathbf{M}(a, t),
\end{aligned}$$

which was to be shown. ∎

**PROOF OF THEOREM 8:** Using the notation from Theorem 7 and its proof (see above), we let $\mathbf{M}_1 = \mathbf{M}(1 - 1/n, -1/n) = \mathbf{C}_{n \times n}$ and $\mathbf{M}_2 = \mathbf{M}(1/n, 1/n) = \mathbf{1}_{n \times n}/n$. Using the weighting functions from the proof of that theorem, we then have:

$$A \equiv A(a,t) = \frac{a(1/n) - (1/n)t}{(1 - 1/n)(1/n) - (1/n)(-1/n)} = a - t = \lambda_*,$$

$$T \equiv T(a,t) = \frac{(1 - 1/n)t - a(-1/n)}{(1 - 1/n)(1/n) - (1/n)(-1/n)} = a - t + nt = \lambda_{**}.$$

The stated result then follows as an application of Theorem 7. ∎